\documentstyle[epsf]{amsart}
\input amssym.def
\input amssym.tex
\newtheorem{thm}{\enspace Theorem}[section]

\def\sminus{\smallsetminus}
\def\conj{
conj}
\def\rpp{\R P^2}

\def\C{\Bbb C}

\def\R{\Bbb R}

\def\abz{\vskip\parskip}

\title{Self-linking number of a real algebraic link}
\author{Oleg Viro}
\address{Department of Mathematics, Uppsala University, S-751 06
Uppsala, Sweden;\abz
Department of Mathematics, University of California,
Riverside, CA 92521;\abz
POMI, Fontanka 27, St.~Petersburg, 191011, Russia}
\email{Viro@@mat.uu.se,  Viro@@math.ucr.edu}
\begin{document}
\maketitle

\begin{abstract} For a nonsingular real algebraic curve in
the 3-dimensional projective space and sphere
a new numeric characteristic is
introduced. It takes integer values, is invariant under rigid isotopy,
multiplied by -1 under mirror reflection. In a sense it is a Vassiliev
invariant of degree 1 and a counter-part of a link diagram writhe.
\end{abstract}

\section{Introduction}\label{sI} In the classical knot theory by a link one
means a smooth closed  1-dimensional submanifold of the 3-dimensional
sphere $S^3$, i.~e. several disjoint circles smoothly embedded into
$S^3$. A classical link may emerge as
the set of real points of a real algebraic curve. First, it gives rise
to questions about relations between invariants of the same curve
which are provided by link theory and algebraic geometry.
Second, it suggests to develop a theory
parallel to the classical link theory, but taking into account the
algebraic nature of objects. From this viewpoint it is more
natural to consider real algebraic links up to isotopy consisting of
real algebraic links, which belong to the same continuous family of
algebraic curves, rather than up to smooth isotopy in the class of
classical links. I call an isotopy of the former kind a {\it rigid
isotopy\/} following a terminology established by Rokhlin \cite{R} in
a similar study of real algebraic plane projective curves and extended
later to various other situations  (see, e.g., \cite{Viro New pr.}).
Of course, there is a forgetting functor: any real
algebraic link can be considered as a classical link and a rigid
isotopy as a smooth isotopy.  It is interesting, how much is lost under
this transition.

In this note I point out a characteristic of a real
algebraic link which is lost.  It is unexpectedly simple. In an obvious
sense it is a nontrivial Vassiliev invariant of degree 1 on the class
of real algebraic knots.\footnote{Recall that a knot is a link
consisting of one component.}  In the classical knot theory the lowest
degree of a nontrivial Vassiliev knot invariant is 2.  Thus there is an
essential difference between classical knot theory and the theory of
real algebraic knots.

The characteristic of real algebraic links which is defined below is
very similar to self-linking number of framed knots. I
call it also {\it self-linking number.\/}  Its definition looks
like a refinement of an elementary definition of the writhe of a knot
diagram, but taking into consideration the imaginary part of the knot.

\section{Self-linking of a nonalgebraic knot}\label{slknonalg}In the
classical theory, a self-linking number of a knot is defined only if the
knot is equipped with an additional structure like framing or just a
vector field nowhere tangent to the knot.\footnote{A framing is a pair
of orthogonal to each other normal vector fields on a knot. There is an
obvious construction which makes a framing from a nontangent vector
field and establishes one to one correspondence between homotopy
classes of framings and nontangent vector fields. The vector fields are
more flexible and relevant to the case.} The self-linking number is the
linking number of the knot oriented somehow and its copy obtained by a
small shift in the direction specified by the vector field. It does not
depend on the choice of orientation, since reversing the orientation of
the knot is compensated by reversing the induced orientation of
its shifted copy. Of course, it depends on the homotopy class of the
vector field.

A knot has no natural preferable homotopy class of framings
which would allow to speak about a self-linking number of the
knot without a special care on choice of a framing.\footnote{Moreover,
the self-linking number is used to define a natural class of framings:
namely, framings with the self-linking number zero.}  Some framings
appear naturally in geometric situations. For example, if one fixes a
generic projection of a knot to a plane, the vector field of directions
of the projection appears. The corresponding self-linking number is
called the {\it writhe\/} of the knot. However, it depends on the
choice of projection and changes under isotopy.

The linking number is  a Vassiliev invariant of order 1 of
two-component oriented links. That means that it changes by a constant
(in fact, by 2) when the link experiences a homotopy with a generic
appearance of an intersection point of the components. Whether the
linking number increases or decreases depends only on the local picture
of orientations near the double point: when it passes from
$\vcenter{\epsffile{f01.eps}}$ through
$\vcenter{\epsffile{f02.eps}}$ to
$\vcenter{\epsffile{f03.eps}}$, the linking number
increases by 2. Generalities on Vassiliev invariants see, e.~g., in
\cite{V}

In a sense the linking number is
the only Vassiliev invariant of degree 1 of two-component oriented
links: any Vassiliev invariant of degree 1 of two-component oriented
links is a linear function of the linking number. Similarly, the
self-linking number is a Vassiliev invariant of degree 1 of framed
knots (it changes by 2 when the knot experiences a homotopy with a
generic appearance of a self-intersection point) and it is the only
Vassiliev of degree 1 of framed knots in the same sense. Necessity of
framing for definition of self-linking number can be formulated now
more rigorously: only constants are Vassiliev invariants of degree 1 of
(nonframed) knots.

The definition of the writhe, which is mimicked below, runs as follows:
for each crossing point of the knot projection one defines a {\it local
writhe\/} equal to $+1$ if near the point the knot diagram looks like
$\vcenter{\epsffile{f03.eps}}$ and $-1$ if it looks like
$\vcenter{\epsffile{f01.eps}}$.
Then one sums up the local writhes over all double points of the
projection. The sum is the writhe.

A continuous change of the projection may cause vanishing of a crossing
point. It happens under the first Reidemeister move shown in
the left hand half of Figure \ref{f1}. This move changes the writhe by
$\pm 1$.

\section{How algebraicity enhances self-linking number}\label{Genrem}
\begin{figure}[h]
\centerline{\epsffile{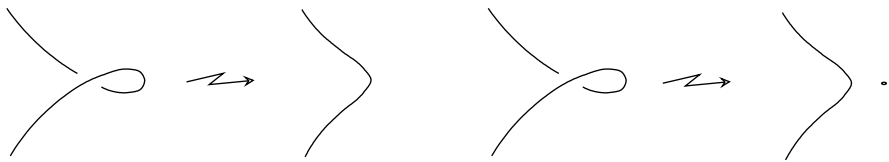}}
\caption{Topological and real algebraic versions of the first
Reidemeister move}
\label{f1}
\end{figure}

If
a link is algebraic then its projection to a plane is algebraic, too. A
generic projection has only ordinary double points.  The total number
of complex double points is constant. The number of real double points
can change, but only by an even number. A real double point cannot turn
alone into imaginary one, as it seems happen under the first
Reidemeister move.  Under the algebraic version of the first
Reidemeister move the double point stays in the real domain, but
becomes solitary, like the only real point of the curve $x^2+y^2=0$.
The algebraic version of the first Reidemeister move is shown in the
right hand half of Figure \ref{f1}.  It is not difficult to prove that
the corresponding family of plane curves can be transformed by a local
diffeomorphism to the family of real rational cubic curves
$y^2=x^2(t-x)$ with $t\in \R$.

A solitary double point of the projection is not image of any real
point of the link. It is the image of two imaginary complex conjugate
points of the complexification of the link. The preimage of the point
in the 3-space under the projection is a real line. It is disjoint from
the real part of the link, but intersects its complexification in a
couple of complex conjugate imaginary points.

In the next section with any solitary double point of the projection,
a local writhe equal to $\pm1$ is associated. It is done in such a way
that the local writhe of the crossing point vanishing in the first
Reidemeister move is equal to the local writhe of the borning solitary
double point.  In the case of an algebraic knot the sum of local
writhes of all double points, both solitary and crossings, does not
depend on the choice of projection and is invariant under rigid
isotopy. This sum is the self-linking number.

There are two types of generic deformations of an algebraic link
changing the rigid isotopy type. One of them is exactly as in the
category of classical links: two pieces of the set of real points come
to each other and pass through each other. A generic projection of
the link experiences an isotopy. No events happen besides that one
crossing point becomes for a moment the image of a double point of the
link and then turns back into a crossing point, but with the opposite
writhe.  Another type has no counter-part in the topological context.
Two complex conjugate imaginary branches pass through each other. At
the moment of passing they intersect in a real isolated double point.
A generic projection of the link experiences an isotopy. No events
happen besides that one solitary double point becomes for a
moment the image of an isolated double point of the link and then turns
back into a usual solitary double point, but with the opposite writhe.

It is clear that the self-linking number of an algebraic knot changes
under both modifications by $\pm2$ with the sign depending only on the
local structure of the modification near the double point. It means
that the self-linking number is the Vassiliev invariant of degree 1.

A construction similar to the construction of the self-linking number
of an algebraic knot can be applied to algebraic {\it links.\/} However
in this case it is necessary either to orient the link or to exclude
from the sum the crossings where the branches belong to distinct
components of the set of real points. In fact, the local writhe depends
on the orientations of the branches, but if the branches belong to the
same component orientations of the branches can be induced from the
same orientation of the component. It is easy to see that the result
does not depend on the choice of orientation of the component.

In the case of knots, self-linking number defines a natural class of
framings, since for knots homotopy classes of framings are enumerated
by their self-linking numbers and we can choose the framing having the
self-linking number equal to the algebraic self-linking number
constructed here.  I do not know any direct construction of this
framing.  Moreover, there seems to be a reason for absence of such a
construction.  In the case of links the construction above gives a
single number, while framings are enumerated by sequences of numbers
with entries corresponding to components.

The construction of this paper can be applied to algebraic links in the
sphere $S^3$. Although from the
viewpoint of knot theory this is the most classical case, from the
viewpoint of algebraic geometry the case of curves in the projective
space is simpler, and I will start from it. The case of spherical
links is postponed to the Section \ref{srS}.

\section{Real algebraic projective links}\label{s0}Let $A$ be a
nonsingular real algebraic curve in the 3-dimensional projective space.
Then the set $\R A$ of its real points is a smooth closed 1-dimensional
submanifold of $\R P^3$, i.~e. a smooth projective link. The set $\C A$
of its complex points is a smooth complex 1-dimensional submanifold of
$\C P^3$.

Let $c$ be a point of $\R P^3$. Consider the projection
$p_c:\C P^3\sminus c\to \C P^2$ from $c$. Assume that $c$ is such that
the restriction to $\C A$ of $p_c$ is generic. This means
that it is an immersion without triple points and at each double point
the images of the branches have distinct tangent lines. As it
follows from well-known theorems, those $c$'s for which this is the
case form an open dense subset of $\R P^3$ (in fact, it is the
complement of a 2-dimensional subvariety).

The real part $p_c(\C A)\cap\R P^2$ of the image consists of the image
$p_c(\R A)$ of the real part and, maybe, several solitary points, which
are double points of $p_c(\C A)$.

There is a purely topological construction which assigns a local writhe
equal to $\pm1$ to a crossing belonging to the image of only one
component of $\R A$. This construction is well-known in the case of
classical knots. Here is its projective version. I borrow it from
Drobotukhina's paper \cite{Dr} on generalization of Kauffman brackets
to links in the projective space.

Let $K$ be a smooth connected one-dimensional submanifold of $\R P^3$,
and $c$ be a point of $\R P^3\sminus K$. Let $x$ be a generic double
point of the projection $p_c(K)\subset \rpp$ and $L\subset \R P^3$ be
the line which is the preimage of $x$ under the projection. Denote by
$a$ and $b$ the points of $L\cap \R P^3$.
\begin{figure}[t]
\centerline{\epsffile{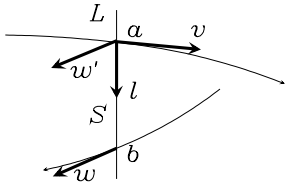}}
\caption{Construction of the frame $v$, $l$, $w'$.}
\label{f2}
\end{figure}
The points $a$ and $b$ divides the line $L$  into two segments.
Choose one of them and denote it by $S$. Choose an orientation of $K$.
Let $v$ and $w$ be tangent vectors of $K$ at $a$ and $b$ respectively
directed along the selected orientation of $K$.
Let $l$ be a vector tangent to $L$ at $a$ and directed inside $S$.
Let $w'$ be a vector at $a$ such that it is tangent to the plane
containing $L$ and $w$ and is directed to the same side of $S$ as $w$
(in an affine part of the plane containing $S$ and $w$). See Figure
\ref{f2}.
The triple $v$, $l$, $w'$ is a base of the tangent space $T_a\R P^3$.
The value taken by the orientation of $\R P^3$ on this frame is the
local writhe of $x$. Its definition involves several choices.  However
it is easy to  prove that the result does not depend on them.

Let $A$, $c$ and $p_c$ be as in the beginning of this Section
and let $s\in\R P^2$ be a solitary double point of $p_c$. Here is a
construction assigning $\pm1$ to $s$. I will call the
result also a {\it local writhe\/} at $s$.

Denote the preimage of $s$ under $p_c$ by $L$. This is a real line in
$\R P^3$ connecting $c$ and $s$.  It intersects $\C A$ in two imaginary
complex conjugate points, say, $a$ and $b$.  Since $a$ and $b$ are
conjugate they belong to different components of $\C L\sminus\R L$.

Choose one of the common points of $\C A$ and $\C L$, say, $a$. The
natural orientation of the component of $\C L\sminus\R L$ defined
by the complex structure of $\C L$ induces orientation on $\R L$ as on
the boundary of its closure. The image under $p_c$ of the local branch
of $\C A$ passing through $a$ intersects the plane of the projection
$\R P^2$ transversally at $s$. Take the local orientation of the plane
of projection such that the local intersection number of the plane and
the image of the branch of $\C A$ is $+1$.

Thus the choice of one of two points of $\C A\cap\C L$ defines an
orientation of $\R L$ and a local orientation of the plane of
projection $\R P^2$ (we can speak only on a local orientation of
$\R P^2$, since the whole $\R P^2$ is not orientable). The plane of
projection intersects\footnote{We may think on the plane of projection
as embedded into $\R P^3$. If you would like to think on it as on the
set of lines of $\R P^3$ passing through $c$, please, identify it in a
natural way with any real projective plane contained in $\R P^3$ and
disjoint from $c$. All such embeddings $\R P^2\to\R P^3$ are isotopic.}
transversally $\R L$ in $s$. The local orientation of the plane,
orientation of $\R L$ and the orientation of the ambient $\R P^3$
determine the intersection number.  This is the local writhe.

It does not depend on the choice of $a$. Indeed, if one  chose
$b$ instead, then both the orientation of $\R L$ and the local
orientation of $\R P^2$ would be reversed. The  orientation of $\R
L$ would be reversed, because $\R L$ receives opposite
orientations from different halves of $\C L\sminus\R L$. The local
orientation of $\R P^2$ would be reversed, because the complex
conjugation involution $\conj:  \C P^2\to\C P^2$ preserves the complex
orientation of $\C P^2$, preserves $\R P^2$ (point-wise) and maps one
of the branches of $p_c(\C A)$ at $s$ to the other reversing its
complex orientation.

Now for any real algebraic projective link $A$ choose a point
$c\in\R P^3$ such that the projection of $A$ from $c$ is
generic and sum up writhes at all crossing points of the projection
belonging to image of only one component of $\R A$ and writhes of all
solitary double points.  The sum is called the {\it self-linking number
of $A$.\/}

It does not depend on the choice of projection. Moreover it is
invariant under {\it rigid isotopy\/} of $A$. By rigid isotopy we
mean an isotopy made of nonsingular real algebraic curves. The effect
of a movement of $c$ on the projection can be achieved by a
rigid isotopy defined by a path in the group of projective
transformations of $\R P^3$.  Therefore the following theorem implies
both independence of the self-linking number on the choice of
projection and its invariance under rigid isotopy.

\begin{thm}\label{mainth} For any two rigidly isotopic real algebraic
projective links $A_1$ and $A_2$ such that their projections from the
same point $c\in\R P^3$ are generic, the self-linking numbers of $A_1$
and $A_2$ defined via $c$ are equal.\end{thm}

To prove this statement, first replace any rigid isotopy by a generic
one. As in purely topological situation of classical links, any generic
rigid isotopy may be decomposed to a composition of rigid isotopies,
each of which makes a local standard move of the projection. There
are 5 local standard moves. They are similar to the
Reidemeister moves. The first of these 5 moves is shown in the right
hand half of Figure \ref{f1}. The next two coincide with the second and
third Reidemeister moves. The fourth move is similar to the second
Reidemeister move: also two double points of projection come to each
other and disappear. However the double points are solitary. The fifth
move is similar to the third Reidemeister move: also
a triple point appears for a moment. But at this triple point only one
branch is real, the other two are imaginary conjugate to each other.
In this move a solitary double point traverses a real branch.

Only in the first, fourth and fifth moves solitary double points are
involved. The invariance under the second and the third move follows
from well-known fact of knot theory that the topological writhe is
invariant under the second and third Reidemeister moves.  Thus we have
to prove that:\begin{enumerate}
\item in the first move the writhe of
vanishing crossing point is equal to the writhe of the borning
solitary point,
\item in the fourth move the writhes of the vanishing
solitary points are opposite and
\item in the fifth move the writhe of the
solitary point does not change.
\end{enumerate}
The proof is not complicated, but would take room inappropriate in this
short note.

The same construction may be applied to real algebraic curves in
$\R P^3$ having singular imaginary points, but no real singularities.
In the construction we can avoid usage of projections from the points
such that some singular point is projected from it to a real point.
Indeed, for any imaginary point there exists only one real line passing
through it (the line connecting the point with its complex
conjugate), thus we have to exclude a finite number of real lines.

\section{Real algebraic links in sphere}\label{srS} The
three-dimensional sphere $S^3$ is a real algebraic variety. It is a
quadric in the four-dimensional real affine space. A stereographic
projection is a birational isomorphism of $S^3$ onto $\R P^3$. It
defines a diffeomorphism between the complement of the center of
projection in $S^3$ and a real affine space.

Given a real algebraic link in $S^3$, one may choose a real point of
$S^3$ from the complement of the link and project the link from this
point to an affine space. Then include the affine space into the
projective space and apply the construction above. The image has no
real singular points, therefore we can use the remark from the end of
the previous section.

\section{Other generalizations}\label{sGeneralizations} It is difficult
to survey all possible generalizations. Here I indicate only two
directions.

 First, consider the most straightforward generalization. Let $L$ be a
nonsingular real algebraic $(2k-1)$-dimensional subvariety in the
projective space of dimension $4k-1$. Its generic projection to $\R
P^{4k-2}$ has only ordinary double points. At each double point either
both branches of image are real or they are imaginary complex
conjugate.  If set of real points is orientable then one can repeat
everything from Section \ref{s0} with obvious changes and obtain a
definition of a numeric invariant generalizing the self-linking number
defined in Section \ref{s0}.

Let $M$ be a nonsingular three-dimensional real algebraic variety with
oriented set of real points equipped with a real algebraic fibration
over a real algebraic surface $F$ with fiber a projective line. There
is a construction which assigns to a real algebraic link (i.~e., a
nonsingular real algebraic curve in $M$) with a generic projection to
$F$ an integer, which is invariant under rigid isotopy, multiplied by
$-1$ under reversing of the orientation of $M$ and is a Vassiliev
invariant of degree 1. This construction is similar to that of Section
\ref{s0}, but uses, instead of projection to $\R P^2$, an algebraic
version of Turaev's shadow descriptions of links \cite{T}.

\end{document}